\documentclass[12pt]{article}
\usepackage[cp1251]{inputenc}
\usepackage{epic}
\usepackage{eepic}
\usepackage{graphics}
\usepackage{amsfonts}
  \usepackage{amsthm}
 \usepackage{amsmath}
 \usepackage{amscd}
 \usepackage{amssymb}
 \usepackage{latexsym}
\usepackage{graphicx}

 \numberwithin{equation}{section}
 \newtheorem{thm}{Theorem}
 
 \newtheorem{prop}{Proposition}
 \theoremstyle{definition}
 
\newcommand{\ds}{\displaystyle}

\title{  Non-Integrability of the Trapped Ionic System}
\author{ Georgi Georgiev \\
 Faculty of Mathematics and Informatics, \\
 Sofia University ``St. Kl. Ohridski'', \\
 5 James Bourchier Blvd., 1164 Sofia, Bulgaria}
\date{}

\begin{document}

\maketitle

\begin{abstract}
In this paper we explore  the two dimensional  system describing trapped
ionic system in the quadrapole field with a superposition of
rationally symmetric hexapole and octopole fields for meromorphic
integrability. We use  the  Lyapunov's and Ziglin-Morales-Ramis
classical methods for the proofs. 
\end{abstract}
\bigskip

\section{Introduction}

Study of the influence of the external fields of the atom occupied a significant place in atomic physics in the early 20th century.
The creation of the capture phenomena by applying static electric and magnetic fields is a remarkable feature of the research of physicists of this period.
The ideal ion trap is based on a pure 3D quadrupole field, on which different types of quadrupole mass spectrometer are based, and the properties of ionic motion are obtained by the exact solution of the resulting Mathieu's equation by analytical methods (see \cite{Trapped} for details).

As it is  known, one of the most used models in non-linear physics is the perturbed harmonic oscillator, because it contains nonlinear behavior that permits testing the different  theories for dynamic systems, as well as its theoretical and experimental applicability in several fields such as particle and plasmas physics (see \cite{Tr13} and \cite{Tr14}), dynamic astronomy -- (\cite{Tr15}, 
\cite{Tr17}) and atomic physics -- (\cite{Tr18}, \cite{Tr19}).
The  mentioned above, the gaps in the experimental configuration and the defects in the physics of the electrodes, lead to the creation of troubles of the multipole field. However, the ion trap is modeled using a two-dimensional oscillator disturbed in harmonic and inharmonic disturbances.

In cylindrical coordinates ($x=r\cos\theta$, $y=r\sin\theta$, $z=z$) and assuming appropriate constraints for simplicity we obtain Hamiltonian
\begin{equation}
\label{1.0}
H=\frac{1}{2}(p_r^2+\frac{p_{\theta}^2}{r^2}+p_z^2)+Ar^2+Bz^2+Cz^3+Dr^2z+Ez^4+Fr^2z^2+Gr^4,
\end{equation}
where $A$, $B$, $C$, $D$, $E$, $F$, and $G$  are  a appropriate
real constants. The Hamiltonian (\ref{1.0}) describes a system with tree degree of freedom having $Z$-axial symmetry -- $\theta$ is
a cyclic and then $p_{\theta}$ is a constant of motion.
The existence of a sufficient number of first integrals of a Hamiltonian system determines whether it is of the two possible types:
quasi-periodic-- integrable or chaotic -- non-integrable. In the case considered in this paper, we have two integrals $H$ and $p_{\theta}$,
and one more is needed for integrability.

We study two dimensional model where $p_{\theta}=0$
\begin{equation}
\label{1.1}
H=\frac{1}{2}(p_r^2+p_z^2)+Ar^2+Bz^2+Cz^3+Dr^2z+Ez^4+Fr^2z^2+Gr^4,
\end{equation}
 for existing an additional  meromorphic integral of motion.

Let  we denote with
$q:=\sqrt{\frac{A}{B}} $, and
$p:=\sqrt{1+\frac{4F}{E}}$. 

The main result of this paper is the following:
\begin{thm}
\label{main2D_1}

a)\, Assume that $q\notin\mathbb{Q} $ or $p\notin\mathbb{Q}$, then 2D system (\ref{1.1})
has no an additional analytic first integral;

b)\, Let $p,\, q\in   \mathbb{Q}$ if $2q\pm p\notin \mathbb{Z}$ or $N(2q)\ge 4$, $N(p)\ge 4$, and $(N(2q),\, N(p))\ne(5,\,5)$  then 2D system (\ref{1.1})
has no an additional meromorphic first integral;

c)\, Let   $2q\pm p\in \mathbb{Z}$, or $N(2q)\le 3$ and  $N(p)\le 3$, or $(N(2q),\, N(p))=(5,\,5)$, 
then 2D system (\ref{1.1}) has no an additional meromorphic first integral if  $D\ne 0$,   ${C}\ne 0$, and $ p^2-1\ne 0$.

 \end{thm}
Here, and in the whole paper, $ N (r) $ is the positive denominator of the irreducible $ r \in \mathbb {Q} $.

The motivation for this work is to look at the problem formulated
in \cite{Trapped}  from another point of view using the
Differential  Galois theory. 

The results of the research  in this work are  almost complete description of the non--integrable cases of the system with Hamiltonian (\ref{1.0}).
The difficulty in the considered system is the presence of seven independent real parameters, which is a great trouble in determining all non--integrable cases.
During the course of the study, the question arises whether the two-dimensional and three-dimensional cases should be considered separately? This distinction is made in \cite{Trapped} , but there the problem is in the difficulty of the calculations in the three-dimensional case.

In the present paper, almost  entire available set of tools of the differential Galois  theory is used to study non-integrability.
We study the variational equation (VE) near a suitable partial  solution for branching (the Lyapunov's method for proving the non-existence of an additional analytical integral), investigation of the commutative  properties of the Galois group 
of VE, and searching for a logarithmic term in the second variational equations.

The  paper s organized as follows:

In section 2 is introduced the two-dimensional model of the considered problem ($p_{\theta}=0$); In sections 3 and susections 3.1  and 3.2, this model is studied for non-integrability;  In section 4, we comment  already obtained  results. There are two
appendices at the end of the text: Appendix A for
 the Ziglin--Morales--Ramis theory and Appendix B for the second order Fuchsian equations and their monodromy.

\section{Statement of the problem}

In this section  we study the two dimensional case ($p_{\theta}=0$) with the Hamiltonian equations:
\begin{align}
\label{1.2}
\dot r  &=  p_r,\, \dot p_r  =  -(2Ar+2Drz+2Frz^2+4Gr^3),\nonumber  \\
\dot z  &=  p_z,\, \dot p_z  =  -(2Bz+3Cz^2+Dr^2+4Ez^3+2Fr^2z).
\end{align}
for existing an additional  integral of motion (here as usual $\dot{}=\frac{d}{dt}$).

In the our study we will suppose that $(F,\, D)\ne (0,\, 0)$, because if we assume the opposite, the variables in the considered system are separated, i.e. the system is integrable.

First, we find a partial solution for (\ref{1.2}). Let we put
$r = p_r = 0$ in (\ref{1.2}) and we have
$$
\ddot z=-(2Bz+3Cz^2+4Ez^3),
$$
multiplying by $\dot z$ and integrating by the time $t$ we have
\begin{equation}
\label{1.3}
{\dot z}^2=-2(Ez^4+Cz^3+Bz^2+h),
\end{equation}
where $h$ is a  constant. Further, we follow the procedures for
Ziglin-Morales-Ramis theory and we find an invariant manifold
$(r,\,p_r,\,z,\, p_z)=(0,\,0,\, z,\, \dot z)$ here $z$ is the
solution of (\ref{1.3}). According to theory, the solution of
(\ref{1.3}) must be a rational function of Weierstrass
$\wp$-function, but it is not important for us right now. Finding
the Variation Equations (VE) we have $\xi_1=dr$, $\eta_1=dp_r$,
$\xi_2=dz$, and  $\eta_2=dp_z$ and we obtain:
\begin{align}
\label{VE1.4}
 \ddot \xi_1 &=  -2(A+Dz+Fz^2)\xi_1,\nonumber\\
\ddot \xi_2  &=  -2(B+3Cz+6Ez^2)\xi_2.
\end{align}
Next we change  the variable in equations (\ref{VE1.4}) by
$\xi_i(t)=\xi_i(z(t))$, where $z(t)$ is a solution of (\ref{1.3}) and
we have $\frac{d\xi_i(t)}{dt}=\frac{d\xi_i}{dz}. \frac{dz(t)}{dt}$,
for $ i=1,\, 2$  and
\begin{eqnarray*}
\frac{d^2\xi_1}{dt^2} &=&\frac{d^2\xi_1}{dz^2}\left (\frac{dz}{dt}\right )^2+\frac{d\xi_1}{dz}\frac{d^2z}{dt^2}\\
&=& -2(Ez^4+Cz^3+Bz^2+h)\frac{d^2\xi_1}{dz^2}-(2Bz+3Cz^2+4Ez^3)\frac{d\xi_1}{dz}\\
&+& 2(A+Dz+Fz^2)\xi_1=0,
\end{eqnarray*}
and
\begin{eqnarray*}
\frac{d^2\xi_2}{dt^2} &=&\frac{d^2\xi_2}{dz^2}\left (\frac{dz}{dt}\right )^2+\frac{d\xi_2}{dz}\frac{d^2z}{dt^2}\\
&=& -2(Ez^4+Cz^3+Bz^2+h)\frac{d^2\xi_2}{dz^2}-(2Bz+3Cz^2+4Ez^3)\frac{d\xi_2}{dz}\\
&+& 2(B+3Cz+6Ez^2)\xi_2=0.
\end{eqnarray*}
If we denote with  $' =\frac{d}{dz}$ we obtain for the VE two
Fuchsian linear differential equations with four singularities:
\begin{eqnarray}
\label{VE1.5}
 \xi_1''  +\frac{4Ez^3+3Cz^2+2Bz}{2(Ez^4+Cz^3+Bz^2+h)}\xi_1'  -\frac{A+Dz+Fz^2}{Ez^4+Cz^3+Bz^2+h}\xi_1=0, \nonumber \\
 \xi_2''  +\frac{4Ez^3+3Cz^2+2Bz}{2(Ez^4+Cz^3+Bz^2+h)}\xi_2'
-\frac{B+3Cz+6Ez^2}{Ez^4+Cz^3+Bz^2+h}\xi_2=0.
\end{eqnarray}
For the Normal Variation Equations (NVE) we suppose that $h=0$ (with
suitable initial conditions for example) and we obtain:
\begin{equation}
\label{NVE}
 \xi_1''  +\frac{4Ez^3+3Cz^2+2Bz}{2(Ez^4+Cz^3+Bz^2)}\xi_1'  -\frac{A+Dz+Fz^2}{Ez^4+Cz^3+Bz^2}\xi_1=0.
 \end{equation}
 The equation ( \ref{NVE}) is a Fuchsian and has four singularities $z=0$, $z=\infty$ and $z=z_i$,
 where $z_{1,2}=\frac{-C\pm\sqrt{C^2-4BE}}{2E}$ ( let $z_1\ne z_2$).  We will drop the index 1 of $\xi_1$.
 \begin{equation}
\label{NVE1}
 \xi''  +\frac{4Ez^2+3Cz+2B}{2z(Ez^2+Cz+B)}\xi'  -\frac{A+Dz+Fz^2}{z^2 (Ez^2+Cz+B)}\xi=0.
 \end{equation}
 The Fuchsian equations with  four or more singularities  like (\ref{NVE1}) are Heun equations.
 The fields of applications of this equations in physics are so large, that it is not possible
to describe them here. However, a examples of many general situations relevant to
physics, chemistry, and engineering where the Heun equations arise  can be found in \cite{Ronveaux}
(pp. 341).

\section{Main result}

Here we study the equation (\ref{NVE1}) for a Liouvillian
solutions. Existence of such solutions of a linear differential equation is equivalent
to the solvability of the identity component of its Galois group.
 
Let us give a simple definition of a Liouville solution of a linear differential equation
\begin{equation}
\label{Lio}
 y''=r(x)y,\, r(x)\in \mathbb{C}[x],
 \end{equation}
 where $\mathbb{C}[x]$ are rational functions with complex coefficients.
 Each second-order linear equation can be written in this way.
 The equation (\ref{Lio}) has Liouville solution, if it is obtained through operations $\int \eta(x)dx$ and $e^{\int \xi (x) dx}$, where $\eta(x)$ and $\xi (x)$ $\in \mathbb{C}[x]$, i.e.
 such a solution is built up by integration and exponentiation.
 
First we find conditions for branching for
solutions of this equation. We will use Lyapunov's classical idea to prove
non-integrability using the branching of the solutions of the
variational equations around an appropriate partial solution. If the
solutions of the VE branching then the Hamiltonian system has not
additional first integral. (See \cite{Lyap} for details.) In this section we will assume that $B\ne 0$ and $E\ne 0$. We find
the  indicial equations to the singular points $z=0$ and $z=\infty$.
We have $\lambda^2-\frac{A}{B}=0$, roots area
$\lambda_j=\pm\sqrt{\frac{A}{B}}$, $j=1,2$ for $z=0$, second one is
$\rho^2-\rho -\frac{F}{E}=0$, with roots
$\rho_k=\frac{1\pm\sqrt{1+\frac{4F}{E}}}{2}$, $k=1,2$ for
$z=\infty$. The solutions around $z=0$ are branching when
$\lambda_j, \, \notin \mathbb{Z}, \, j=1,\, 2 $ and near $z=\infty$
we have branching - when $\rho_k, \, \notin \mathbb{Z}, \, k=1,\,2
$. This is an application of the Frobenius's method (\cite{Poole} p.
70): the fundamental system of solutions around of the singular
points are presented in the form
 $\Phi_j(z)=z^{\lambda_j}\Psi_j(z)$, where $\Psi_j(z)$ are holomorphic functions (locally) for $j=1,2$. Then the solutions of (\ref{NVE1})
 are branching for $\lambda_j, \, \notin \mathbb{Z}, \, j=1,\, 2 $. For $z=\infty$   we do the same and obtain that $\rho_k, \, \notin \mathbb{Z}, \, k=1,\,2 $. It should be noted that it is possible, if $\lambda_j\in \mathbb{Q}$ and $\rho_k\in \mathbb{Q}$, we also can achieve not branching by a standard changing of the variables. This change of variables does not affect the commutativity of the unity component of the Galois group.
We have proved the following proposition.
\begin{prop}
\label{Lyap}
Let $q=\sqrt{\frac{A}{B}}\notin\mathbb{Q} $ or $p=\sqrt{1+4\frac{F}{E}}\notin\mathbb{Q}$, then 2D system (\ref{1.2})
has no additional holomorphic first integral.
\end{prop}

 We have completed the first step. Now we use the notation
$F=\ds{\frac{p^2-1}{4}}E$, $A=q^2B$, where
$p,\, q \in \mathbb{Q}$. (This is the case when the solutions are not branching i. e.  $\lambda_j, \, \rho_k \in \mathbb{Q} $.)

 In eq. (\ref{NVE1}) we put 
 $$a(z):=\frac{4Ez^2+3Cz+2B}{2z(Ez^2+Cz+B)}=\frac{1}{z}+\frac{\frac{1}{2}}{z-z_1}+\frac{\frac{1}{2}}{z-z_2},$$
$$b(z):= -\frac{Fz^2+Dz+A}{z^2 (Ez^2+Cz+B)}=\frac{-q^2}{z^2}+\frac{\alpha}{z}+\frac{\beta}{z-z_1}+\frac{\gamma}{z-z_2},$$
where
$$\alpha = \frac{D}{B}+\frac{AC}{B^2},$$
$$\beta =-\frac{Fz_{1}^2+Dz_1+A}{Ez_{1}^2(z_1-z_2)},$$
and
$$\gamma =\frac{Fz_{2}^2+Dz_2+A}{Ez_{2}^2(z_1-z_2)}.$$
In the terminology of \cite{Chur} we get 
$a_{\infty}=\lim_{z\to\infty} z a(z)=2$, $b_{\infty}=\lim_{z\to\infty} z^2 b(z)=\ds{\frac{1-p^2}{4}}$, $\Delta_{\infty}=\sqrt{(1-a_{\infty})^2-4b_{\infty}}=\pm p$, and $a_0=1$, $b_0=-q^2$, $\Delta_0=\sqrt{(1-1)^2-4(-q^2)}=\pm 2q$. Then, we get $t_0=-2\cos (\pi\Delta_0)=-2\cos( \pi 2q)$, and $t_{\infty}=2\cos (\pi\Delta_{\infty})=2\cos (\pi p)$. The values of $\Delta_1$ and $\Delta_2$ corresponding to points $z_1$ and $z_2$ are much more complex as expressions and depend on too many parameters.

The first result of our exposition is the following proposition:
\begin{prop}
\label{2Dmain}
Let  $q=\sqrt{\frac{A}{B}}$ and $p=\sqrt{1+\frac{4F}{E}}$ are rational numbers, then 2D
system (\ref{1.2}) has no an additional meromorphic first integral when $2q\pm p\notin \mathbb{Z}$ and $N(2q)\ge 4$, $N(p)\ge 4$, and $(N(2q),\, N(p))\ne(5,\,5)$.

(Here $N(r)$ is the positive denominator of the irreducible fraction $r\in \mathbb{Q}$.)
\end{prop}

{\bf Proof:}
The assumptions in the Proposition are equivalent to condition $\mathbb{Q}[\cos\pi 2q]\ne \mathbb{Q}[\cos\pi p]$ (this is a consequence of the Theorem \ref{Berger}). We obtain that $t_{\infty}\notin \mathbb{Q}[t_0]\subset \mathbb{Q}[t_0,\, t_1,\, t_2]$, i. e. $t_{\infty}$ is transcendental over $\mathbb{Q}[t_0,\, t_1, \, t_2]$.
  The Proposition \ref{2Dmain} is a direct application of Theorem \ref{BaChurch1} and  Theorem \ref{Berger}.

  This finishes the proof of Proposition \ref{2Dmain}.

\subsection{$VE_2$ }

In this section we assume that the  monodromy group is abelian, i. e.
  $2q\pm p\in \mathbb{Z}$, or $N(2q)\le 3$, $N(p)\le 3$, or $(N(2q),\, N(p))=(5,\,5)$, and we proceed to the second variation,
which will give us additional conditions for non-integrability. The procedure is standard, we have:
 \begin{eqnarray*}
r &= & \varepsilon\xi_{11}+\varepsilon^2 \xi_{21}+\dots , \\
z & = & z(t)+\varepsilon\xi_{12}+\varepsilon^2 \xi_{22}+\dots , \\
 p_r & = & \varepsilon\eta_{11}+\varepsilon^2 \eta_{21}+\dots , \\
 p_z & = & \dot{z}(t)+\varepsilon\eta_{12}+\varepsilon^2 \eta_{22}+\dots ,
 \end{eqnarray*}
here $(p_r,\, r,\, p_z,\,  z)=(0,\, 0,\, \dot{z}(t), \,z(t))$ is an invariant manifold of the system (\ref{1.2}).
We substitute in the system  changing the variables $t\rightarrow z(t)$ and we obtain
\begin{eqnarray}
\label{VE1.5_1}
 \xi_{11}''  +\frac{4Ez^2+3Cz+2B}{2z(Ez^2+Cz+B)}\xi_{11}'  -\frac{A+Dz+Fz^2}{z^2(Ez^2+Cz+B)}\xi_{11}=0, \nonumber \\
 \xi_{12}''  +\frac{4Ez^2+3Cz+2B}{2z(Ez^2+Cz+B)}\xi_{12}'
-\frac{B+3Cz+6Ez^2}{z^2(Ez^2+Cz+B)}\xi_{12}=0,
\end{eqnarray}
\begin{eqnarray}
\label{VE2}
 \xi_{21}''  & + & \frac{4Ez^2+3Cz+2B}{2z(Ez^2+Cz+B)}\xi_{21}'  -\frac{A+Dz+Fz^2}{z^2(Ez^2+Cz+B)}\xi_{11} =K_2^{(1)}\nonumber\\
 \xi_{22}''  & + & \frac{4Ez^2+3Cz+2B}{2z(Ez^2+Cz+B)}\xi_{22}'
-\frac{B+3Cz+6Ez^2}{z^2(Ez^2+Cz+B)}\xi_{22} = K_2^{(2)}.
\end{eqnarray}
We have
$$
K_2^{(1)}=\frac{2Fz+D}{z^2(Ez^2+Cz+B)}\xi_{11}\xi_{12}=\tilde{K}_2^{(1)}\xi_{11}\xi_{12},
$$
$$K_2^{(2)}=\frac{2Fz+D}{z^2(Ez^2+Cz+B)}\xi_{11}^2+\frac{12Ez+3C}{z^2(Ez^2+Cz+B)}\xi_{12}^2,$$
and
 $$f_2=(0,\, K_2^{(1)},\, 0,\, K_2^{(2)})^T.$$
We put in (\ref{VE1.5_1}) and (\ref{VE2})
$F=\ds{\frac{p^2-1}{4}}E$, $A=q^2B$, $2q\pm p\in \mathbb{Z}$, or $N(2q)\le 3$, $N(p)\le 3$, or $(N(2q),\, N(p))=(5,\,5)$.
First we  find  linearly independent solutions close to $z_i$, ($i=1,\, 2$) of (\ref{VE1.5_1}) - $\xi_{11}^{(1)}$, $\xi_{11}^{(2)}$ and $\xi_{12}^{(1)}$, $\xi_{12}^{(2)}$.
Without losing a community we can assume that
$\xi_{11}^{(1)}\dot{\xi_{11}^{(2)}}-\xi_{11}^{(2)}\dot{\xi_{11}^{(1)}}=1$ and $\xi_{12}^{(1)}\dot{\xi_{12}^{(2)}}-\xi_{12}^{(2)}\dot{\xi_{12}^{(1)}}=1$ .
Then the fundamental matrix of (\ref{VE1.5_1}) and its inverse are
\begin{equation}
\label{X}
X (z) =
\begin{pmatrix}
 \xi_{11} ^{(1)} & \xi_{11} ^{(2)} & 0 & 0\\
 \dot{\xi}_{11} ^{(1)} & \dot{\xi}_{11} ^{(2)}   & 0 & 0 \\
0 & 0 & \xi_{12} ^{(1)} & \xi_{12} ^{(2)}\\
0 & 0 & \dot{\xi}_{12} ^{(1)} & \dot{\xi}_{12} ^{(2)}
  \end{pmatrix},
  \end{equation}

 \begin{equation}
\label{X^1}
X ^{-1}(z) =
\begin{pmatrix}
 \dot{\xi}_{11} ^{(2)} & -\xi_{11} ^{(2)} & 0 & 0\\
 -\dot{\xi}_{11} ^{(1)} & \xi_{11} ^{(1)}   & 0 & 0 \\
0 & 0 & \dot{\xi}_{12} ^{(2)} & -\xi_{12} ^{(2)}\\
0 & 0 & -\dot{\xi}_{12} ^{(1)} & \xi_{12} ^{(2)}
  \end{pmatrix}.
  \end{equation}
  We will show that  a logarithmic term
appears in  local solution of (${{VE}}_2$). For this purpose,
it is sufficient to show that  at least one component of $X^{-1} f_2$
has a nonzero residue at $z_1$. We calculate
 of $X^{-1} f_2$, which looks like
 $$(-\xi_{11}^{(2)}K_2^{(1)},\,\xi_{11}^{(1)}K_2^{(1)},\, -\xi_{12}^{(2)}K_2^{(2)},\, \xi_{12}^{(1)}K_2^{(2)} )^T.$$
Now we find proper solutions of (\ref{VE1.5_1})  in the neighbourhood of $z=z_1$ and we have:
 \begin{equation}
\label{xi_11^1}
\xi_{11}^{(1)}=c_1(z-z_1)^{1/2}\left(1+\frac{(p^2-4)z_1+(4q^2+2)z_2+4\frac{D}{E}}{6z_1(z_1-z_2)}(z-z_1)+\dots \right),
  \end{equation}
 \begin{equation}
\label{xi_11^2}
\xi_{11}^{(2)}=c_2\left(1+\frac{(p^2-1)z_1+4q^2z_2+4\frac{D}{E}}{2z_1(z_1-z_2)}(z-z_1)+\dots \right),
  \end{equation}
 \begin{equation}
\label{xi_12^1}
\xi_{12}^{(1)}=l_1 (z-z_1)^{1/2} \left(1+\frac{3z_1-2z_2}{2z_1(z_1-z_2)}(z-z_1)+\dots \right),
  \end{equation}
\begin{equation}
\label{xi_12^2}
\xi_{12}^{(2)}=l_2  \left(1+\frac{6z_1-4z_2}{z_1(z_1-z_2)}(z-z_1)+\dots \right),
  \end{equation}

   \begin{equation}
\label{K_2^1}
\tilde{K}_2^{(1)}
 =  \left(\frac{\frac{p^2-1}{4}z_1+\frac{D}{E}}{z_1^2(z_1-z_2)}(z-z_1)^{-1}+\dots \right).
  \end{equation}
 We  choose the constants $c_1$, $c_2$, $l_1$ and $l_2$ so that the Wronskians of these pairs of solutions are a units.

 {\bf Remark.} We obtain that the Galois group of the second part of the equation (\ref{VE1.5_1}) is represented by the matrix group
 $\left\{\begin{pmatrix}
1 & 0  \\
\mu & 1
\end{pmatrix},  \,\mu\ne0\right\}$ which is
commutative.

 We  write expressions for residue  of
 $\tilde{K}_2^{(1)}\xi_{11}^{(1)}\xi_{12}^{(1)}\xi_{11}^{(2)}$  at $z=z_1$:
  \begin{equation*}
\label{Res1}
Res_{z=z_1}(\tilde{K}_2^{(1)}\xi_{11}^{(2)}\xi_{12}^{(2)}\xi_{11}^{(2)} )
 =  \frac{\frac{p^2-1}{4}z_1+\frac{D}{E}}{z_1^2(z_1-z_2)},
  \end{equation*}
and  for $z=z_2$ we have
\begin{equation*}
\label{Res2}
  Res_{z=z_2}(\tilde{K}_2^{(1)}\xi_{11}^{(2)}\xi_{12}^{(2)}\xi_{11}^{(2)} )
 =  \frac{\frac{p^2-1}{4}z_2+\frac{D}{E}}{z_2^2(z_1-z_2)}.
\end{equation*}

The conditions $\frac{p^2-1}{4}z_i+\frac{D}{E}=0 $, $i=1,\, 2$ are fulfilled when $D= 0$, $C=0$ and $p^2-1= 0$ (for $z_1\ne z_2$).
 
  We proved the following:
  \begin{prop}
\label{mainVE2}
Let   $2q\pm p\in \mathbb{Z}$, or $N(2q)\le 3$ and  $N(p)\le 3$, or $(N(2q),\, N(p))=(5,\,5)$, 
then 2D system (\ref{1.2}) has no an additional meromorphic first integral if is true that 
$ D\ne 0$, $C\ne 0$, and $p^2-1\ne 0$.
\end{prop}
This complete the study of the two-dimensional case.
\subsection{The Degenerate Cases}

Here we consider some cases in the study of equations in variations (NVE)  we omitted so far  to simplify the presentation. The cases we need to look at in this section are $z_1=z_2$, $A=B=0$, $F=E=0$, $B=0$ and $A\ne 0$, $E=0$, and $F\ne 0$.

During the study of equation (\ref{NVE}), we assumed that $z_1\ne z_2$,  therefore we  need to consider the case $z_1= z_2$.
Then we receive $C^2=4BE$ (the discriminant) and
  \begin{equation}
\label{NVE1.1}
 \xi''  +\frac{2(C^2z^2+3CBz+2B^2)}{C^2z(z-\frac{2B}{C})^2}\xi'  -\frac{4(A+Dz+Fz^2)}{C^2z^2 (z-\frac{2B}{C})^2}\xi=0.
 \end{equation}
 This is the Heun's confluent equation, which we transform to standard form ($\ds{\frac{d^2Y}{dx^2}}=r(x)Y$) using the M\"{o}bius transformation
 $z=\ds{\frac{-x+1}{\frac{C}{2B}x}}$ and we have
 \begin{equation}
 \label{CHE}
 \frac{d^2Y}{dx^2}=\left (9-\frac{\frac{2D}{C}}{x}-\frac{\frac{2D}{C}}{1-x}+\frac{\frac{F}{E}}{x^2}+\frac{\frac{A}{B}+\frac{1}{4}}{(1-x)^2}\right )Y.
 \end{equation}
 The issue of the conditions under which Liouville's solutions of such equations exist is studied in detail in \cite{BaBi}.
 Using the notation: we have $\alpha^2=36$, $\eta=\ds{\frac{1}{2}-\frac{2D}{C}}$, $\delta=0$, $\beta^2=1+\ds{\frac{4F}{E}}$,
 and $\gamma^2=1+\ds{\frac{4A}{B}}$.
 The conditions that equation (\ref{CHE}) has no Liouville solutions are $\pm\beta\pm\gamma\notin\left(\mathbb{Z}_{even}\setminus \{0\}\right )
 $ ($\alpha\ne 0$ and $\delta=0$), which  expressed by the parameters of our problem are
 $\pm\sqrt{1+\ds{\frac{4F}{E}}}\pm \sqrt{1+\ds{\frac{4A}{B}}}\notin (\mathbb{Z}_{even}\setminus \{0\})$.
 We prove the following theorem
\begin{prop}
\label{det0}
When $C^2=4BE$,  the system with Hamiltonian (\ref{1.1}) is not meromorphic integrable  iff
 $$\pm\sqrt{1+\ds{\frac{4F}{E}}}\pm \sqrt{1+\ds{\frac{4A}{B}}}\notin (\mathbb{Z}_{even}\setminus \{0\}).$$
\end{prop}

 For the last case  which we have to consider is when in  (\ref{NVE}) we put $C=0$ and $E=0$.  In (\ref{NVE})  we get
   \begin{equation}
\label{NVE1.2}
 \xi''  +\frac{1}{z}\xi'  -\frac{4(A+Dz+Fz^2)}{Bz^2}\xi=0,
 \end{equation}
 and we obtain  a standard form for (\ref{NVE1.2})
 \begin{equation}
 \label{Witt1}
 \frac{d^2Y}{dx^2}=\left (\frac{F}{B}+\frac{\frac{4A}{B}+1}{4x^2}+\frac{\frac{D}{B}}{x}\right )Y.
 \end{equation}
 After replacing  $x=\ds{\frac{1}{2}\sqrt{\frac{B}{F}}.T}$ we obtain  a Whittaker equation
  \begin{equation}
 \label{Witt2}
 \frac{d^2Y}{dT^2}=\left (\frac{1}{4}+\frac{\frac{A}{B}+1}{4T^2}+\frac{\frac{D}{2B}\sqrt{\frac{B}{F}}}{T}\right )Y.
 \end{equation}
In the notation of \cite{MR2} we have
$$\kappa=-\frac{D}{2B}\sqrt{\frac{B}{F}},$$
and
$$\mu=\pm\sqrt{\frac{1}{2}+\frac{A}{B}}.$$
The conditions for non-integrability of the Whittaker equation are
 $$\pm\kappa\pm\mu=\pm\frac{D}{2B}\sqrt{\frac{B}{F}}\pm\sqrt{\frac{1}{2}+\frac{A}{B}}\notin (\mathbb{Z}+\frac{1}{2}).$$
 (See \cite{MR2} for details.)
We will prove
\begin{prop}
\label{det0.1}
When $C=0$ and $E=0$,  the system with Hamiltonian (\ref{1.1}) is not meromorphic integrable  iff
 $$\pm\frac{D}{2B}\sqrt{\frac{B}{F}}\pm\sqrt{\frac{1}{2}+\frac{A}{B}}\notin (\mathbb{Z}+\frac{1}{2}).$$
\end{prop}
Conditions $\pm\sqrt{1+\ds{\frac{4F}{E}}}\pm \sqrt{1+\ds{\frac{4A}{B}}}\in (\mathbb{Z}_{even}\setminus \{0\})$ for $C^2=4BE$,  and $\pm\frac{D}{2B}\sqrt{\frac{B}{F}}\pm\sqrt{\frac{1}{2}+\frac{A}{B}}\in (\mathbb{Z}+\frac{1}{2})$ for $C=E=0$ can also be studied with a second variations, which will remain an open question.

We continue with the case $A=B=F=E=0$, which is integrable (the solutions of this system can be written explicitly). 

Let us now $A\ne 0$, $B\ne0$, $D\ne 0$ and $E=F=0$, then the  Hamilton function and the equations of the system are:
\begin{equation}
\label{HamF=E=0}
H=\frac{1}{2}(p_r^2+p_z^2)+Ar^2+Bz^2+Cz^3+Dr^2z+Gr^4
\end{equation}
\begin{align}
\label{SysF=E=0}
\dot r  &=  p_r,\, \dot p_r  =  -(2Ar+2Drz+4Gr^3),\nonumber  \\
\dot z  &=  p_z,\, \dot p_z  =  -(2Bz+3Cz^2+Dr^2).
\end{align}
For $p_r=r=0$ we get a partial solution $z=z(t)$ that satisfies the differential equations
$$\ddot{z}=-z(3Cz+2B)$$ 
and
$${\dot{z}}^2=-2Cz^3-2Bz^2+h.$$
The last equation has  a solution $z(t)=\ds{\frac{\wp(t)+B/6}{-C/2}}$, here $g_2=\ds{\frac{4B^2}{9}}$, $g_3=\ds{-\frac{B^3}{18}-\frac{C^2}{4}h}$, and
$$g_2^3-27g_3^2=\frac{13}{2916}B^6-\frac{3}{4}B^3C^2H-\frac{27}{16}C^4h^2\ne0.$$

For the first and second variations, we obtain sequentially
\begin{eqnarray*}
 \ddot\xi_{11}   +  2(A+Dz)\xi_{11}=0, \nonumber \\
 \ddot\xi_{12}   + 2(B+3Cz)\xi_{12}=0,
\end{eqnarray*}
\begin{eqnarray*}
\ddot\xi_{21}   +   2(A+Dz)\xi_{21} =-D\xi_{11}\xi_{12}\nonumber\\
\ddot\xi_{22}   +  2(B+3Cz)\xi_{22} =-3C{\xi_{12}}^2 -D{\xi_{11}}^2.
\end{eqnarray*}
The  equation for $\xi_{11}$ is of the Lam\'e type. Let's write it in the form
$$ \ddot\xi_{11}   -  (4\frac{D}{C}\wp(t)+\frac{BD}{C}-2A)\xi_{11}=0.$$
We make the following notation $\frac{4D}{C}=n(n+1)$ and we obtain
$$ \ddot\xi_{11}   -  (n(n+1)\wp(t)+B\frac{n(n+1)}{4}-2A)\xi_{11}=0.$$

We use the notation in Appendix B1 and we get $\alpha(t,\,h):=n(n+1)\wp(t)+B\frac{n(n+1)}{4}-2A$,
\begin{equation*}
P(\alpha ,\, h):= (a_1 + h a_2) \alpha^3 + (b_1 + h b_2) \alpha^2 + (c_1 + h c_2) \alpha + (d_1 + h d_2),
\end{equation*}
and we obtain the coefficients 
\begin{equation*}
a_1=\frac{4}{n(n+1)},\,\,a_2=0,
\end{equation*}
\begin{equation*}
b_1=-3B+\frac{6A}{n(n+1)},\,\, b_2=0,
\end{equation*}
\begin{equation*}
c_1=\frac{11}{36}n(n+1)B^2+\frac{3A^2}{n(n+1)}-3AB,\,\, c_2=0,
\end{equation*}
\begin{equation*}
d_1=\frac{A^3}{2n(n+1)}+\frac{5n^2(n+1)^2}{48}B^3-\frac{3}{4}A^2B+\frac{11n(n+1)}{72}B^2A,\,\, d_2=\frac{C^2}{4}.
\end{equation*}
The conditions for non-integrability in this case are the negations of the conditions in Theorem \ref{thA} for the coefficients $a_j$, $b_j$, $c_j$, and $d_j$, $j=1,\,2$. 
We consider a second variation for the case $D=3C$ (where condition 1 of Theorem \ref{thA} is satisfied), i.e. for $n=3$ (or $n=-4$).
In the neighborhood of $t=0$ we obtain the following solutions of VE1:
\begin{equation*}
\xi_{11}^{(1)}(t):=t^{4}\left(1-\frac{(B-A)}{9}t^2+\left(\frac{A^2}{198}-\frac{AB}{99}+\frac{B^2}{90}\right)t^2+\dots \right),
  \end{equation*}
 \begin{equation*}
\xi_{11}^{(2)}(t):=\frac{1}{t^3}+\frac{(A-B)}{5}\frac{1}{t}+\frac{(3A^2-6BA+B^2)}{90}t+\dots ,
  \end{equation*}
 \begin{equation*}
\xi_{12}^{(1)}(t):= t^4 \left(1+\frac{B}{9}t^2+\frac{B^2}{90}t^4+\dots\right),
  \end{equation*}
\begin{equation*}
\xi_{12}^{(2)}(t):= \frac{1}{t^3}-\frac{B}{5}\frac{1}{t}+\frac{B^2}{90}t+\dots,
  \end{equation*}

Then the residue in $t=0$ for $D\xi_{11}^{(2)}\xi_{12}^{(2)}\xi_{11}^{(2)}$ is 
\begin{equation*}
D\left(\frac{A^4}{900}-\frac{8}{1125}A^3B+\frac{97}{6750}A^2B^2-\frac{34}{3375}AB^3+\frac{23}{13500}B^4\right). 
 \end{equation*}
 The above expression is equal to 0 just when $A=B$ (because $D\ne 0$ ).
Then  the condition for a non-zero logarithmic term is $A\ne B$. We proved the following proposition:
\begin{prop}
\label{det0.11}
When  $E=F=0$, $D=3C$ ($n=3$ or $n=-4$),  the system with Hamiltonian (\ref{1.1}) is not meromorphic integrable  if
 $A\ne B$.
\end{prop}
{\bf Remark:} Unfortunately, there are too many cases that can be investigated with second variations in the assumptions of the above Proposition. But they can be formulated as open questions. 

The study of case $ A = B = 0 $ is similarly to already considered general case, and we will only summarize its results below
\begin{prop}
\label{det0.12}
When  $A=B=0$ and $D=\ds{\frac{P^2-1}{4}}C$, $F=\ds{\frac{p^2-1}{4}}E$  then the system with Hamiltonian (\ref{1.1}) is not  integrable  if

a) at least one of $P$ or $p$ is not a rational number;

b) for $P$, $p$ $\in \mathbb{Q}$ one of the sub-cases is fulfilled:

b1) $P\pm p \notin\mathbb{Z} $;

b21) $P\pm p \in\mathbb{Z}$ and $N(P)\ge 4$ and $N(p)\ge 4$;

b22) $P\pm p \in\mathbb{Z}$ and $N(P)= 5$ and $N(p)= 5$;

b3) $P\pm p \in\mathbb{Z}$ and $N(P)\le 3$ and $N(p)\le 3$ and $E\ne 0$ and  
$$3P^2p^2-P^4-6p^2-P^2+5\ne 0.$$

\end{prop}
We can prove the following
\begin{prop}
\label{det0.12}
When  $A=0$ and $B\ne 0$ and $2FC\ne DE$ then the system with Hamiltonian (\ref{1.1}) is not meromorphic integrable.
\end{prop}
\begin{prop}
\label{det0.13}
When  $E=0$ and $F\ne 0$,  then the system with Hamiltonian (\ref{1.1}) is not meromorphic integrable if 
$2F\ne C$ or $16F\ne 3C$.
\end{prop}

Another different approach in studying integrability of the two-dimensional problem is possible, it can be seen in \cite{DarbouxPoint}. The problem is reduced to a Hamiltonian system with homogeneous potentials, and has a well-developed research scheme. Let us explain what this method is about and how we can use it for limit  number of integrable cases.
We consider a Hamiltonian system 
 \begin{equation}
 \label{HomPot}
H=\frac{1}{2}\left(p_1^2+p_2^2\right)+V(r,\,z) ,
  \end{equation}
with potential 
  \begin{equation*}
  V(r,\,z)=V_{min}(r,\,z)+\dots+ V_{max}(r,\,z),
  \end{equation*}
  which is a sum of homogeneous potentials. Here with $V_{min}(r,\,z)$ (respectively $V_{max}(r,\,z)$ ) we mean the smallest ( largest ) possible degree of homogeneous part of $V(r,\,z)$. Potential $V(r,\,z)$ is called integrable if its corresponding Hamiltonian system (\ref{HomPot}) is integrable. As it is noted in \cite{Jarno}, if $V(r,\,z)$ is an integrable potential, then $V_{min}(r,\,z)$ and $V_{max}(r,\,z)$ are also integrable. In our case we have $V_{min}(r,\,z)=Ar^2+Bz^2$ and $V_{max}(r,\,z)=Ez^4+Fr^2z^2+Gr^4$ and since $V_{min}(r,\,z)$ is integrable, then the possible integrable cases are those for which $V_{max}(r,\,z)$ is integrable. These potentials are fully investigated in \cite{DarbouxPoint} and using the notation in this paper the possible integrable cases in $V_{max}(r,\,z)$ are: $V_1$, $V_3$, $V_4$, $V_5$ and $V_6$. Let us consider these cases in the context of the fourth degree potential in our task. 
  
 $V_1:$
 
 a) $F=E=0$ we have considered this case above, and unfortunately there have been too many possible integrable cases ( the negations of Theorem \ref{thA} ) for the coefficients found there.
 
 b) $F=G=0$ this case is integrable (possibly) for $\sqrt{\frac{A}{B}}\in\mathbb{Q}$ and $D=0$ (because $p^2=1$).

 $V_3:$  In this case we have $G=E=\frac{1}{4}$, $F=\frac{1}{2}$ here the necessary conditions for integrability are $\sqrt{\frac{A}{B}}\in\mathbb{Q}$, $C=0$ and $D=0$.
 
  $V_4:$  We have $E=1$, $F=0$, condition for integrability is only $\sqrt{\frac{A}{B}}\in\mathbb{Q}$ and $D=0$ (here $p^2=1$).
  
  $V_5:$
  
  a) $E=\frac{1}{4}$, $F=3$, $G=4$
  
  b) $E=4$, $F=3$, $G=\frac{1}{4}$ in both cases the necessary conditions for integrability are $\sqrt{\frac{A}{B}}\in\mathbb{Q}$, $C=0$ and $D=0$.
  
    $V_6:$
    
    a) $E=\frac{1}{4}$, $F=\frac{3}{2}$, $G=2$
  
  b) $E=2$, $F=\frac{3}{2}$, $G=\frac{1}{4}$ in both cases the necessary conditions for integrability are $\sqrt{\frac{A}{B}}\in\mathbb{Q}$, $C=0$ and $D=0$.
    
 The obtained necessary conditions for integrability can be studied numerically, this can be done in a separate study.
 
This completes our study on most degenerative cases.

\section{Remarks and Comments}
Let us now consider the  results in paper \cite{Trapped} (for 2D- case) in the context of already obtained conclusions.
It should be noted that the additional first integrals (page 4 in \cite{Trapped}) found in this article are not fully true  (see \cite{Comment}).

In the case noted there with 1(i) we have $q^2=1$, $p=\pm 5$ and $2q\pm p\in \mathbb{Z}$ and because $D=3C$ in my opinion this case should not be integrable if $C\ne 0$ or $E\ne 0$. For $C=0$  the integrability is possible. The case when $E=F=0$ leads us to the result of 2b(i).
 In case 1(ii-b) (in notations of  \cite{Trapped}) with additional condition $A=B$ according to \cite{Comment}, we have $q^2=1$, $p=\pm 1$ and $2q\pm p\in \mathbb{Z}$  therefore the condition with the second variations is not fulfilled and then we have an integrability and additional integral. In the case 2b(i) we have the  conditions from Proposition \ref{det0.11}  and we have an  integrability  with an additional integral (with small corrections, see \cite{Comment}). In the first of  2b(ii) -- we have $q^2=4$, $p=\pm 5$ and 
$2q\pm p\in \mathbb{Z}$ and this is $V_6$-case - the condition with the second variations is not fulfilled--$D=3C=0$ and then we may have integrability. In the last of 2b(ii) we have $q^2=4$, $p=\pm 7$ and 
 $2q\pm p\in \mathbb{Z}$ and this is $V_5$-case-- here we can expect integrability, since the conditions for the second variations are not fulfilled  $D=3C=0$.
 
 We could assume  there are  more cases in which integrability is possible (eventually) such as when $2q\pm p\in \mathbb{Z}$, or $N(2q)\le 3$, $N(p)\le 3$, or $(N(2q),\, N(p))=(5,\,5)$, $D=0$, and $p^2-1=0$ or $C=0$. In fact, this case is exactly when $F=E=D=0$ and this is the  case left unconsidered by us that is integrable because, the variables are separated.

There is a number of unexplored cases (for example for $F=E=0$, if the conditions of Theorem \ref{thA} or the negations of Propositions \ref{det0}, \ref{det0.1}, \ref{det0.12} and \ref{det0.13} are satisfied) that, although atypical (degenerated), could present a complete picture of the chaos in this two-dimensional problem. Unfortunately, the study of each such case is too complicated as a job due to the large number of free parameters (7) in the initial model.
 
 In conclusion, we can say that the research in \cite{Trapped} gives an interesting and qualitative result. The inaccuracies admitted by the authors are something normal when performing such a large volume of analytical and symbolic calculations, which in my opinion do not reduce the quality of their result.

 \section*{ Acknowledgments }

The author has been partially supported by RMS N
577/17.08.2018/pr.11 from Bulgarian Government and grant  80-10-43/
26.03.2020 from Sofia University ``St. Kliment Ohridski''.

\vspace{5ex}

\noindent
{\Large \bf APPENDICES}

\appendix

 \section{Ziglin--Morales--Ramis Theory}

In this section we  recall some classical and more advanced results that we used for  our research.

We say that a system (with $n$ degree of freedom) is integrable in the sense of Liouville, if it has a complete set of  $n$ independent
first integrals in involution. 
 We recall some notions and facts about integrability of Hamiltonian systems in the complex domain,
the Ziglin--Morales--Ramis theory and its relations with differential Galois groups of linear equations. We will follow \cite{MR1} and \cite{ChG}.

We consider a Hamiltonian system
\begin{equation}
\label{2.1}
\dot{x} = X_{H} (x), \quad t \in \mathbb{C}, \quad x \in M
\end{equation}
corresponding to an analytic Hamiltonian $H$, defined on the complex
$2 n$-dimensional manifold $M$. If we suppose the system (\ref{2.1}) has a
non-equilibrium solution $\Psi (t)$, we denote by $\Gamma$ its phase
curve. We can write the equation in variation (VE) near this solution
\begin{equation}
\label{2.2}
\dot{\mathbf{\xi}} = D X_{H} ( \Psi (t)) \mathbf{\xi},
\quad \mathbf{\xi} \in T_{\Gamma} M.
\end{equation}

Further, we consider the normal bundle of $\Gamma$, $F:= T_{\Gamma} M / TM$ and let
$\pi : T_{\Gamma} M \to F$ be the natural projection. The equation (\ref{2.2})
leads to an equation on $F$
\begin{equation}
\label{2.3}
\dot{\eta} = \pi_{*} (D X_{H} ( \Psi (t))(\pi^{-1} \eta) , \quad \eta \in F.
\end{equation}
which is  called a normal variational equation (NVE) around $\Gamma$.
The (NVE) (\ref{2.3}) recognizes a first integral $d H$, linear on the fibers of $F$.
The level set $F_r := \{\eta \in F | d H (\eta) = r \} , r \in \mathbb{C}$,
 is $(2 n - 2)$-dimensional affine bundle over $\Gamma$.
We will call $F_r$ the reduced phase space of (\ref{2.3}) and the restriction of
the (NVE) on $F_r$ is called the reduced normal variational equation.

Then the main result of the Morales--Ramis \cite{MR1} theory is:
\begin{thm}
\label{th2}
Let's assume that the Hamiltonian system (\ref{2.1}) has $n$ meromorphic first integrals in involution,
then the identity component $G^0$ of the Galois group of the variational equation is abelian.
\end{thm}

Next we consider a linear  system
\begin{equation}
\label{2.10}
y' = A (x) y, \quad y \in \mathbb{C}^n ,
\end{equation}
or  linear homogeneous differential equation, which is essentially the same
\begin{equation}
\label{2.11}
y^{(n)} + a_1 (x) y^{(n-1)} + \ldots + a_n (x) y = 0,
\end{equation}
with $x \in \mathbb{CP}^1$  and $A \in \mathrm{gl} (n, \mathbb{C} (x))$,
($a_j (x) \in \mathbb{C} (x))$.
Let $S:=\{x_1, \ldots, x_s\}$ be the set of singular points of (\ref{2.10}) (or (\ref{2.11})) and let $Y (x)$ be a fundamental solution
of (\ref{2.10})  (or (\ref{2.11}))  at $x_0 \in \mathbb{C} \setminus S$. Acording to existence theorem  this solution is analytic
near of $x_0$. The continuation of $Y (x)$
 along a nontrivial loop on $\mathbb{CP}^1$ defines  a linear automorphism
of the space of the solutions, called the monodromy. Analytically
this transformation can be presented as a follows: the linear automorphism $\Delta_{\gamma}$,
associated with a loop $\gamma \in \pi_1 (\mathbb{CP}^1 \setminus S, x_0)$ corresponds
to multiplication of $Y (x)$ from the right by a constant matrix $M_{\gamma}$, called monodromy matrix
$$
\Delta_{\gamma} Y (x) = Y (x) M_{\gamma}.
$$
The set of these matrices forms the monodromy group.

We add another object to the  (\ref{2.10}) (or (\ref{2.11})) - a differential
Galois group.
We have a differential field $K$, that is a field with a derivation $\partial = '$, i.e.
an additive mapping satisfying derivation rule. Differential automorphism
of $K$ is an auto\-mor\-phism commuting with the derivation.

The coefficient field in (\ref{2.10}) (and  (\ref{2.11})) is $K = \mathbb{C} (x)$. Let $y_{i j}$ be
elements of the fundamental matrix $Y (x)$. Let $L (y_{i j})$ be the
extension of $K$ generated by $K$ and $y_{i j}$ -- a differential
field. This extension is called a Picard--Vessiot's extension.
 Similarly to classical Galois Theory we define the Galois group
$G := Gal_{K} (L) = Gal (L/K)$ to be the group of all differential
automorphisms of $L$ leaving the elements of $K$ fixed.
 Galois group is an algebraic group. It has an unique
connected component $G^0$ which contains the identity and  is a normal
subgroup of finite index.  Galois group $G$ can be represented
as an algebraic linear subgroup of $\mathrm{GL} (n, \mathbb{C})$ by
$$
\sigma (Y (x)) = Y (x) R_{\sigma},
$$
where $\sigma \in G$ and $R_{\sigma} \in \mathrm{GL} (n, \mathbb{C})$.

We can do the same locally at $a \in \mathbb{CP}^1$, replacing $\mathbb{C} (x)$ by the field of germs of meromorphic
functions at $a$. In this way we can speak of a local differential Galois group $G_a$ of (\ref{2.10}) at $a \in \mathbb{CP}^1$,
defined in the same way for Picard-Vessiot extensions of the field $\mathbb{C} \{x-a\}[(x-a)^{-1}]$.

It should be noted that by its definition the monodromy group is contained
in the differential Galois group of the corresponding system.

Next, we present some facts from the theory of linear systems with singularities.
We call a singular point $x_i$  regular if any of the solutions of (\ref{2.10})
(or of (\ref{2.11})) has at most polynomial growth in arbitrary sector with a vertex at
$x_i$. Otherwise the singular point is called  irregular.

We say that the system (\ref{2.10}) has a singularity of the Fuchs type at $x_i$ if $A (x)$ has a
simple pole at $x = x_i$. For the equation (\ref{2.11}) the Fuchs type singularity at  $x_i$
means that the functions $(x - x_i)^j a_j (x)$ are holomorphic in a neighborhood of $x_i$.

If the system (\ref{2.10}) has a singularity of the Fuchs type, then  this singularity is regular.
The opposite is not true. However, for the equation (\ref{2.11}) the regular singularities coincide with
the singularities of  Fuchs type.

A system with only regular singularities  is called Fuchsian system.
For such systems we have :
\begin{thm}
\label{th3}
({\bf Schlesinger} )
The differential Galois group coincides with the  Zariski closure in $\mathrm{GL}(n,\mathbb{C})$ of the mo\-no\-dro\-my
group.
\end{thm}

The fact that $G^0$ is abelian  doesn't imply  necessarily
integrability of the Hamiltonian system.
There is a method which, in the case of abelian Galois group, can draw conclusion when the system (\ref{2.1}) is non-integrable. This method based on
the higher variational equations has been introduced in \cite{MR1}
and  the Theorem \ref{th2} has been extended in \cite{MRS1}. What is the idea of higher variational
equations? For the system (\ref{2.2}) with a particular solution
$\Psi (t)$ we put
\begin{equation}
\label{2.5}
x = \Psi (t) + \varepsilon \xi^{(1)} + \varepsilon^2
\xi^{(2)} + \ldots + \varepsilon^k \xi^{(k)} + \ldots,
\end{equation}
where $\varepsilon$ is a small   parameter. When substituting the
above expression into Eq. (\ref{2.2}) and comparing terms with the
same order in $\varepsilon$ we obtain the following chain of
linear non-homogeneous equations
\begin{equation}
\label{2.6}
\dot{\xi}^{(k)} = A (t) \xi^{(k)} + f_k (\xi^{(1)},
\ldots, \xi^{(k-1)}), \quad k = 1, 2, \ldots ,
\end{equation}
where $A (t) = D X_{H} ( \Psi (t))$ and $f_1 \equiv 0$. The
equation (\ref{2.6}) is called k-th variational equation
(${\rm{VE}}_k$). Let $X (t)$ be the fundamental matrix of
(${\rm{VE}}_1$)
$$
\dot{X} = A (t) X .
$$
Then the solutions of $({\rm{VE}}_k), k > 1$ can be found by
\begin{equation}
\label{2.7}
\xi^{(k)} = X (t) c (t),
\end{equation}
where $c (t)$ is a solution of
\begin{equation}
\label{2.8}
\dot{c} = X^{-1} (t) f_k .
\end{equation}
Although (${\rm{VE}}_k$) are not actually homogeneous equations,
they can be placed in this framework, and therefore,  successive
extensions $ K \subset L_1 \subset L_2 \subset \ldots \subset
L_k$ can be defined, where $L_k$ is the extension obtained by adjoining the
solutions of (${\rm{VE}}_k$). The
differential Galois groups $Gal (L_1 /K), \ldots, Gal (L_k /K)$ can be defined accordingly. The following
result is proven in \cite{MRS1}.
\begin{thm}
({\bf  Morales-Ruiz, Ramis,  Sim\'o})
\label{th3}
If the Hamiltonian system (\ref{2.2}) is integrable
in Liouville sense then the identity component of each Galois
group $ Gal (L_k /K)$ is abelian.
\end{thm}

Note that we apply Theorem \ref{th3} to the situation where the
identity component of the Galois group $Gal(L_1/K)$ is abelian.
This means that the first variational equation is solvable. Once
we have the solution of $({\rm{VE}}_1)$, then the solutions of
$({\rm{VE}}_k)$ can be found by the method of 
constant variations as explained above. Hence, the differential Galois groups $Gal(L_k /K)$ are solvable.
One possible way to show that some of them are not commutative is
to find a logarithmic term in the corresponding solution.
We need to explain why the existence of a non-zero logarithmic term in $VE_k$ around some singular point guarantees us non-integrability. The Galois group $Gal(L_k /K)$ is abelian, if and only if,  the local monodromy of the $({\rm{VE}}_k)$ around the singular point of the coefficients is identity. If for some $k$ ($k=2$ in our case), we obtain non-zero residue in the Laurent expansions of the expressions of $ X^{-1} (t) f_k$, near singularity point, then the local monodromy will be represented by  a lower (or upper) triangular matrix which is not an identity, i. e. the Galois group $Gal(L_k /K)$ is not abelian (see
detailed descriptions and explanations in \cite{MR1,MRS1,MR2}).

\section{ Fuchsian Differential Equations of the Second Order}

In \cite{Chur}  a  result is proved for the monodromy group of second order Fuchsian differential equations  with rational coefficients. Let us try to explain this in a little more details. We consider the equation
\begin{equation}
\label{Heun2}
 \frac{d^2 Y(x)}{dx^2}  + C_1(x)\frac{d Y(x)}{dx}
   + C_2(x)Y(x)=0,
 \end{equation}
where
$$C_k=\sum_{j=1}^m\frac{\lambda_j}{(x-\alpha_j)^{l_k^j}}+\dots,\, k=1,\,2.$$
(Here ${l_k^j}$ is the order of the pole $\alpha_j$.) We use the following notations $A_j:=\lambda_j^1(C_1)$, $B_j:=\lambda_j^2(C_2)$ for $j=1,\dots m$, $A_{\infty}:=\lim_{x\to\infty}(x-\alpha_j)^{l_k^1} C_1(x)$ and $B_{\infty}:=\lim_{x\to\infty}(x-\alpha_j)^{l_k^2} C_2(x)$, we use also notation  $\Delta_j:=\sqrt{(A_{j}-1)^2-4B_{j}}$, $\Delta_{\infty}:=\sqrt{(A_{\infty}-1)^2-4B_{\infty}}$ and $t_j:=-2\cos\pi \Delta_j$, and $t_{\infty}:=2\cos\pi \Delta_{\infty}$ then we have next 
\begin{thm}
({\bf Baider--Churchill})
\label{BaChurch1}
Suppose that $t_{\infty}$ is transcendental over $\mathbb{Q}[t_1,\, t_2,\dots \, t_m]$, then the monodromy group is not abelian.
\end{thm}
In fact, \cite{Chur} shows a possible connection between Classical and Differential Galois theory.

This result has been known since 1991, but its use was difficult because checking for transcendence is not a simple job. The study of the real cyclotomic extensions of the field from rational numbers turns out to be the difficulty here, and more precisely the significant  question  of the transcendence of $\cos \pi r_1$ and $\cos \pi r_2$ over $\mathbb{Q}$.
Fortunately, in 2018 the issue was resolved in \cite{Berger} in the form of the following
\begin{thm}
({\bf Berger})
\label{Berger}
Let $r_1, \, r_2\in \mathbb{Q}$ be such that neither $r_1\pm r_2$ is an integer. Then the following are equivalent:

1) The numbers 1, $\cos \pi r_1$ and $\cos \pi r_2$ are $\mathbb{Q}$ independent;

2)  $N(r_j)\ge 4$ for $j=1,\, 2$, and $(N(r_1),\, N(r_2))\ne(5,\,5)$.
\end{thm}

This statement allows us to determine when two real cyclotomic fields coincide.

\subsection{Necessary conditions for integrability of Hamiltonian systems which have (NVE) of Lam\'e type}

We recall some facts about the integrability of Hamiltonian systems with two degrees
of freedom, an invariant surface and which (NVE) are of Lam\'e type. We follow the exposition in \cite{MSim,MR1}.

The Lam\'e equation is written in the form
\begin{equation}
\label{Lame}
\ddot{\xi} - (n(n+1) \wp(t) + B) \xi = 0,
\end{equation}
where $\wp (t)$ is the Weierstrass function with $g_2$ and $g_3$, satisfying
$\dot{v}^2 = 4 v^3 - g_2 v - g_3$ with $\Delta = g_2 ^3 - 27 g_3 ^2 \neq 0$.

The known  cases of closed form solutions of (\ref{Lame}) are:

(i) The Lam\'e and Hermite solutions. In this case $n \in \mathbb{Z}$ and $g_2, g_3, B$ are arbitrary parameters;

(ii) The Brioschi-Halphen-Crowford solutions. Here $m:= n + 1/2 \in \mathbb{N}$ and the parameters
$g_2, g_3, B$ must satisfy an algebraic equation.

(iii) The Baldassarri solutions. Now
$n + 1/2 \in \frac{1}{3} \mathbb{Z} \cup \frac{1}{4} \mathbb{Z} \cup \frac{1}{5} \mathbb{Z} \setminus \mathbb{Z}$ with
additional algebraic relations between the other parameters.

Note that in the case (i) the identity component of the Galois group $G^0$ is of the form
$\begin{pmatrix}
1 & 0 \\
\nu & 1
\end{pmatrix}$
 and in the cases (ii) and (iii) $G^0 = id$ ($G$ is finite).
And these are the all cases when the Lam\'e equation is integrable.

Now consider  two degrees of freedom Hamiltonian
\begin{equation}
\label{ham2}
H = \frac{1}{2} (p_1 ^2 + p_2 ^2) + V (q_1, q_2),
\end{equation}
$q_j (t) \in \mathbb{C}, p_j (t) = \dot{q}_j, j=1, 2$. We assume that  a
family of solutions  exists there
$$
\Gamma_h : q_2 = p_2 = 0, \quad q_1 = q_1 (t, h), \quad  p_1 (t, h) = \dot{q}_1 (t, h)
$$
and $q_1 (t, h)$ is a solution of
$$
\frac{1}{2} \dot{q}_1 ^2 + \varphi (q_1) = h, \quad h \in \mathbb{R} .
$$
The (NVE) around $\Gamma_h$ is
\begin{equation}
\label{nve}
\ddot{\xi} - \alpha (t, h) \xi = 0 ,
\end{equation}
where $\alpha (t, h) = \alpha (q_1 (t, h))$ is such that (\ref{nve}) is of type (\ref{Lame}).

In \cite{MSim,MR1} the type of the potentials $V$ with this property are obtained as well as
the necessary conditions for the integrability of the Hamiltonian systems with the Hamiltonian
(\ref{ham2}). In order to formulate the result we need  additional notations.

Since $\alpha (t, h)$ depends linearly on $\wp (t)$, then $\dot{\alpha}^2$ is a cubic polynomial in
$\alpha$, depending also in $h$
\begin{equation}
\label{apl}
\dot{\alpha}^2 : = P (\alpha, h) = P_1 (\alpha) + h P_2 (\alpha) .
\end{equation}
The following coefficients are introduced
\begin{equation}
\label{coeff}
P (\alpha, h) = (a_1 + h a_2) \alpha^3 + (b_1 + h b_2) \alpha^2 + (c_1 + h c_2) \alpha + (d_1 + h d_2) .
\end{equation}
Note that the following Theorem gives necessary conditions only
from the analysis of the first variational equation.
\begin{thm}
\label{thA}
({\bf Morales-Ruiz-Simo}). Assume that a natural Hamiltonian system has (NVE) of Lam\'e type,
associated to the family of solutions $\Gamma_h$, lying on the plane $q_2 = 0$ and parametrized by the
energy $h$. Then, a necessary conditions for integrability is that the related polynomials $P_1$ and
$P_2$ satisfy $a_2 = 0$, and one of the following conditions holds:

\vspace{2ex}

\noindent
1. $a_1 = \frac{4}{n (n+1)}$ for some $n \in \mathbb{N}$;

\vspace{2ex}

\noindent
2. $a_1 = \frac{16}{4 m^2 - 1}$ for some $m \in \mathbb{N}$. Then, if assumption the conjecture above is true, one should
have $b_2 = 0$ and we should be in one of the following cases:

\vspace{1ex}

2.1) $m = 1$ and $b_1 = 0$,

\vspace{1ex}

2.2) $m = 2$ and $c_2 = 0, \, 16 a_1 c_1 + 3 b_1 ^2 = 0$,

\vspace{1ex}

2.3) $m = 3$ and $16 a_1 d_2 + 11 b_1 c_2 = 0, \, 1024 a_1 ^2 d_1 + 704 a_1 b_1 c_1 + 45 b_1 ^3 = 0$,

\vspace{1ex}

2.m) $m > 3$. Then, we should have $b_1 = 0$ and, furthermore, either $c_1 = c_2 = 0$ if $m$ is
congruent with $1, 2, 4$ or $5$ modulo $6$, or $d_1 = d_2 = 0$ if $m$ is odd;

\vspace{2ex}

\noindent
3.  $a_1 = \frac{4}{n (n+1)}$ with
$n + 1/2 \in \frac{1}{3} \mathbb{Z} \cup \frac{1}{4} \mathbb{Z} \cup \frac{1}{5} \mathbb{Z} \setminus \mathbb{Z}$,
$b_2 = 0$ and either $c_2 = 0, b_1 ^2 - 3 a_1 c_1 = 0$
or $c_2 b_1 - 3 a_1 d_2 = 0, 2 b_1 ^3  - 9 a_1 b_1 c_1 + 27 a_1 ^2 d_1 = 0$.
\end{thm}
It is easy to notice that the condition 1. in the above Theorem gives the Lam\'e and Hermite solutions (i), the condition 2.--
the Brioschi-Halphen-Crowford solutions (ii), and the condition 3. -- the Baldassarri solutions (iii).

\end{document}